\documentclass[a4paper,11pt]{article}

\usepackage[dvips]{graphicx} 
\usepackage[latin1]{inputenc}
\usepackage{amssymb}         
\usepackage{amsmath}         
\usepackage{amsfonts}
\usepackage[usenames,dvipsnames]{xcolor}
\usepackage{curves}
\usepackage{a4wide}

\newcommand{\nn}{\nonumber}

\begin{document}

\title{String Stability towards Leader thanks to Asymmetric Bidirectional Controller}


\author{Arash Farnam\thanks{SYSTeMS research group, Faculty of Engineering and Architecture, Ghent University; Technologiepark Zwijnaarde 914, 9052 Zwijnaarde(Ghent), Belgium (e-mail: arash.farnam@ugent.be).}
 and Alain Sarlette\thanks{SYSTeMS research group and  project team, INRIA Paris, France (e-mail: alain.sarlette@inria.fr)}} 

\maketitle

\begin{abstract}                
This paper deals with the problem of string stability of interconnected systems with double-integrator open loop dynamics (e.g.~acceleration-controlled vehicles). We analyze an asymmetric bidirectional linear controller, where each vehicle is coupled solely to its immediate predecessor and to its immediate follower with different gains in these two directions. We show that in this setting, unlike with unidirectional or symmetric bidirectional controllers, string stability can be recovered when disturbances act only on a small ($N$-independent) set of leading vehicles. This improves existing results from the literature with this assumption. We also indicate that string stability with respect to arbitrarily distributed disturbances cannot be achieved with this controller. 
\end{abstract}




\section{Introduction}\label{sec:intro}

The platooning problem is both a practical and theoretical topic for automated vehicles, that includes many different issues. One of the benchmark settings, commonly called the \emph{vehicle chain}, is relevant e.g.~for automated highway systems, see e.g.~\cite{1,2,3,5,Jova,new}.
In this setting, a set of vehicles are arranged on a single path and their objective is to keep a desired distance with respect to their predecessor and follower, while the first vehicle additionally has to track a commanded trajectory. The main issue is the behavior of the chain when the number of vehicles $N$ becomes very large. The open-loop model of each vehicle is a double integrator, in accordance with positions as outputs and forces $\simeq$ accelerations as input. Most of the numerous methods to design distributed controllers for this interconnected system can guarantee input-to-output stability, but would still lead to increasingly big oscillations of the vehicle chain with increasing number of vehicles, which has been formalized among others as string instability.

Since its definition in (\cite{5}; \cite{6}), string (in)stability has attracted a lot of discussion. Recently researchers have characterized a lot of details and variants on the issue, to the point that this conference paper can only offer a truncated view of the literature. The following papers are just the closest ones to the problem at hand, and we must apologize for leaving out probably tens of significant papers which are just farther from our focus. Essentially, it has been established as an unavoidable shortcoming of linear controllers that none of them can guarantee string stability, in several precisely identified distributed control settings.

In the simplest setting, when vehicles look at relative velocities and relative positions of their preceding vehicle, the transfer function from vehicle $i$ to $i+1$ takes the form of a complementary sensitivity function. It then follows from the Bode integral that any stable linear controller always leads to a transfer function with an $\infty$-norm more than one, and thus an exponential growth of an initial disturbance of some frequency as it travels along the vehicle chain (\cite{15}; \cite{5}; \cite{6}). As the chain grows longer, the last vehicle (with index $N$) will thus undergo larger and larger oscillations. The absence of an $N$-independent bound on these oscillations is what we here call $L_2$ string instability, and it further implies what we here call $(L_2,l_2)$ string instability, namely the sum of squares of the motions of all the vehicles is unbounded. The $(L_2,l_2)$ string instability becomes important when small disturbances can act on all the vehicles: if a disturbance input on a single vehicle implied a bounded yet non-vanishing effect on the whole chain, then when small disturbances act on all the vehicles this effect would sum up to become unbounded on each vehicle as $N$ grows.

The above fundamental result has been extended to the case where each vehicle looks at a limited number of `neighbor' vehicles \emph{in front} of them (\cite{1}; \cite{2}; \cite{3}; \cite{5}; \cite{6}). Another line of work has considered \emph{bidirectional} coupling --- i.e.~each vehicle can react to some vehicles just in front and to some vehicles just behind itself. When the coupling is symmetric, the mutual reactions of two interconnected vehicles can be modeled following mechanical principles --- e.g.~placing a suitably tuned spring-damper system between them, and analyzing it with passivity type methods. It has been shown with this approach that the impact of a bounded input disturbance on the error in the distance between any single pair of vehicles can be kept bounded with a suitable design (\cite{17}), i.e.~$L_2$ string stability can be achieved. Yet $(L_2,l_2)$ is impossible to achieve, i.e.~for any linear symmetric bidirectional controller looking only one vehicle in front and one vehicle behind, the $l_2$ norm of the vector of distance errors will necessarily grow unbounded for some $l_2$-bounded input disturbances on the vehicles (\cite{15}; \cite{Barnew}).

The present paper is concerned with asymmetric bidirectional coupling, where the vehicle reacts differently to its predecessor than to its follower in the chain. The benefit, on a different objective, of breaking the symmetry in the coupling has been famously shown in \cite{7}. Unfortunately, some limitations of this setting have also been proved. If the asymmetry just consists of a constant factor in front of the controller \cite{TAC-Martinec}, then $l_2$ string stability will fail. Furthermore, it has been established that a PD controller cannot work and that keeping symmetric DC controller gain is a necessary condition for string stability \cite{Automatica-Knorn}.

Our contribution rather follows up on the more positive observations in \cite{MartinecEJC}. Like in this paper, we consider an asymmetric PD coupling where \emph{disturbances act on the first vehicle(s)} only. Our analysis also turns out to follow a similar flow-inspired analysis. We add two more positive properties to the setting of  \cite{MartinecEJC}, namely:\newline
(i) the system with these assumptions satisfies not only $L_2$ but also $(L_2,l_2)$ string stability\newline
(ii) more detailed analysis shows that there is no need to worry about flow reflections at the end of the chain, so no need to introduce a dedicated controller on the last vehicle.\newline
While these observations do not solve the practical problem of $(L_2,l_2)$ string instability when disturbances can act on \emph{any} vehicle, they might form a valuable basis when minimal variations on the setting are sought towards achieving this goal. 

The impossibility results discussed above hold for vehicles modeled as second-order pure integrators and relying on purely relative measurements. We probably must mention that a successful line of work has shown how adding a term proportional to absolute velocity to the dynamics, can solve the string instability problem. In proposed solutions, this absolute velocity can take the form of a drag force or introduced in the actual controller, e.g.~in what has become known as the time headway policy or adaptive cruise control (\cite{Dirk}; \cite{new}; \cite{10}; \cite{16}; \cite{new}). We believe that despite these results, the theoretical interest in achieving string stability without absolute velocity remains justified for practical purposes. In some applications at least (e.g.~space flight, underwater), one might question the availability of a reliable, globally accessible common reference with respect to which the absolute velocity of all the vehicles can be measured. Moreover, relying on drag to ensure a positive property is probably not the best control engineering solution, when modern transportation systems like the latest vacuum tube transit proposal (see e.g.\cite{19}) try to minimize the drag for energy efficiency purposes. 

The paper is organized as follows. Section \ref{2} presents the setting. Section \ref{3} contains its detailed analysis and the main result, while Section \ref{4} illustrates it with simulations. 

\noindent \textbf{Acknowledgment:} The authors have to thank an anonymous reviewer for sharing their very clear viewpoint on recent string stability investigations.


\section{Model description}\label{2}

\subsection{String stability, general}

The $H_{\infty}$ norm of transfer function $C(s)$ is given by $\Vert C \Vert_{\infty}=sup_{\omega\ge 0}\vert C(j\omega)\Vert$. \emph{Re} and \emph{Im} respectively denote the real and imaginary parts. 

Consider a family $\{\mathcal{S}_N \}_{N=1,2,...}$ of networks. Each network $\mathcal{S}_N$ consists of $N+1$ interconnected dynamical subsystems, whose configuration we denote by $x(t) = (x_0(t),x_1(t),x_2(t),...,x_N(t))$ and which can be subject to input disturbances $\, d(t) = (d_0(t),d_1(t),d_2(t),...,d_N(t))$. The focus of this work lies on the \emph{relative} states of the subsystems with respect to each other, while their absolute value remains free. More precisely, we assume that the coordinates have been chosen such that the control objective is to stabilize the subspace $x_0=x_1=x_2=...=x_N$. The actual value of $x_0$ can then be independently guided as e.g.~a trajectory tracking command. The context of vehicle chains considers the most basic network topology, where subsystem $k$ is coupled to the subsystems $k-1$ and $k+1$, for $k=1,2,3,...,N-1$. The formal objective of string stability reflects this topology in the configuration error vector $e(t) = (e_1(t),e_2(t),...,e_N(t))$ with each $e_k = x_{k-1} - x_k$. There are several variants of string stability in the literature, and as explained in the introduction we here go for the stronger one. The $(L_2, l_2)$ norm of a time-dependent vector e.g.~$x(t)$ is defined by
$$\Vert x(\cdot) \Vert_2 = \left( \; \sum_{k=0}^N  \int_{-\infty}^{+\infty} \vert x_k(t) \vert^2 dt \; \right)^{1/2} \, .$$
\vspace{2mm}

\noindent \textbf{Definition}: \emph{The family of networks $\{\mathcal{S}_N \}_{N=1,2,...}$ is $(L_2, l_2)$ string stable if for every $\epsilon > 0$ there exists $\delta > 0$ such that: $\Vert d(\cdot) \Vert_2 < \delta$ implies $\Vert e(\cdot) \Vert_2 < \epsilon$ for all networks i.e.~all $N=1,2,...\;$ .}\vspace{2mm}

In other words, the focus of string stability is that the configuration error must be bounded \emph{uniformly in $N$}. The weaker notion of $L_2$ string stability requests a uniform bound for all $e_k$, instead of taking the sum over subsystems. A fully realistic comparison however is between $\Vert d(\cdot) \Vert_2 < N\, \delta$ and $\Vert e(\cdot) \Vert_2 < N\, \epsilon$ when disturbances can affect any vehicle. Then as the sum goes both over the subsystems and over time, it is not enough for string stability to e.g.~evacuate an input disturbance by transporting it towards the tail of the chain: in addition, the disturbance must be damped at a rate that is bounded away from zero. The $(L_2,l_2)$ criterion furthermore allows a standard analysis in frequency domain, through Parseval's equality, involving e.g.~$H_\infty$ norms of transfer functions.

\subsection{Vehicle chain}

String stability has been the focus of major interest in the following model by \cite{5}. Consider $N$ vehicles modeled as pure double-integrators:
\begin{eqnarray}\label{system}
\ddot{x}_k(t)=u_k+d_k \   ,    \    k=0,1,2,...,N \, .
\end{eqnarray} 
Here $x_k$ is the absolute position of vehicle $k$, while $u_k$ and $d_k$ are acceleration control input and disturbance input, respectively. The objective of each vehicle is to follow its preceding vehicle at a fixed desired distance $r$, in appropriate coordinates $x_k \longrightarrow x_k-kr$ this can be reformulated as stabilizing $x_0=x_1= ... = x_N$. To achieve this task, vehicle $k$ adapts $u_k$ as a function of observed information about its neighboring vehicles. We introduce two fundamental assumptions about this information.
\begin{itemize}
\item[(A1)] The feedback controller $u_k$ can only depend on \emph{relative} states of the vehicles, e.g.~their relative positions  $x_k-x_{k-1}$ or relative velocities  $\dot{x}_k-\dot{x}_{k-1}$.
\item[(A2)] The controller $u_k$ of vehicle $k$ can only depend on such information from a few neighboring vehicles, i.e.~whose index is comprised in $[k-\bar{k},k+\bar{k}]$ for some (small) $\bar{k}$ independent of $N$.
\end{itemize}
Furthermore, we impose that the controller of a given vehicle $k$ should not depend on $N$. This means in essence that the vehicle applies its control action only by looking at its local neighborhood, without knowing anything about the rest of the chain (except tacitly acknowledging that they will all cooperate).
It is under these assumptions that fundamental impossibilities to obtain string stability with linear controllers have been established, as explained in the introduction.

\begin{figure}[h!]\label{fig1}
\includegraphics[width=0.5\textwidth, trim=165mm 0mm 140mm 297mm, clip=true]{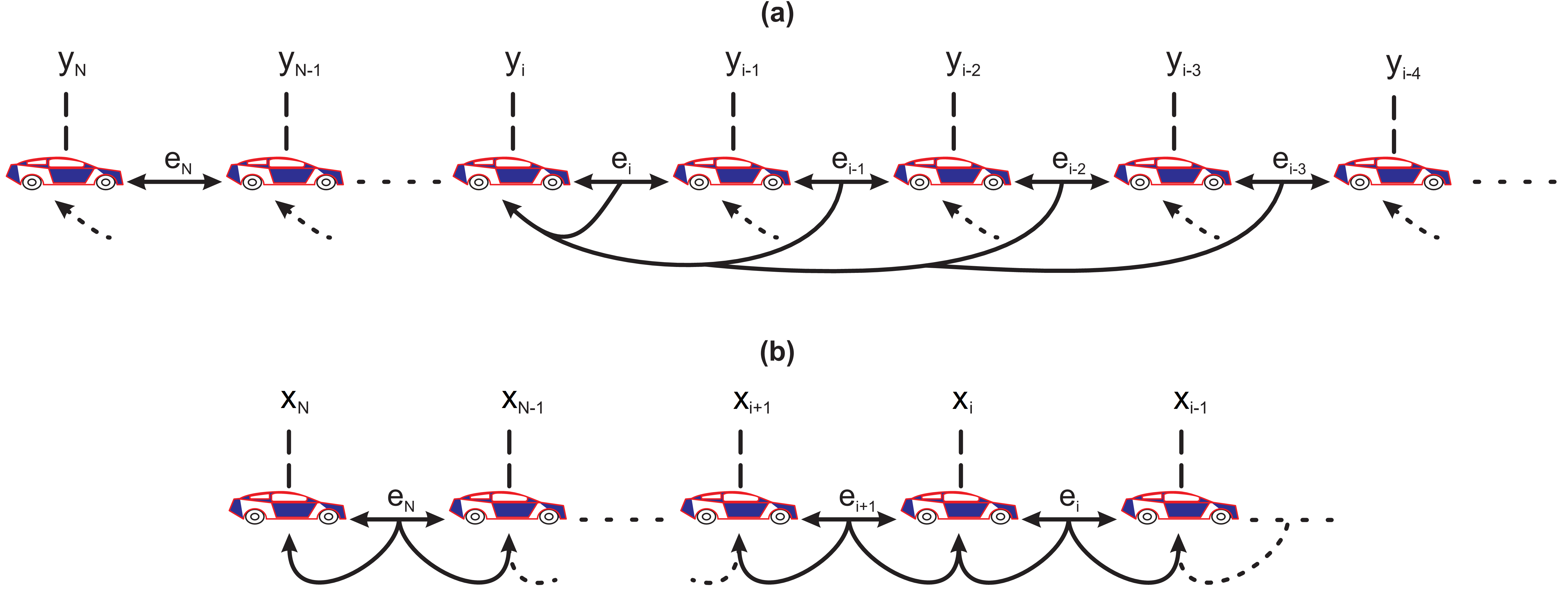}  
\caption{Vehicle chain with bidirectional coupling to closest neighbors.}\label{1}
\end{figure}

The present paper considers the model \eqref{system} with assumptions (A1) and (A2), more precisely vehicle $k$  relies on relative information about one preceding vehicle $k-1$ and one following vehicle $k+1$. The scheme of this controller is shown on Fig.~\ref{1}. Like in \cite{7,TAC-Martinec,Automatica-Knorn,MartinecEJC}, the feedback transfer function assigned to the preceding vehicle can differ from the feedback transfer function assigned to the following vehicle (asymmetry), and the point of our paper is to highlight the benefits of this asymmetry.

\subsection{A simple asymmetric controller}\label{3}

Explicitly, we consider the control:
\begin{eqnarray}\label{eq2}
u_0 & = &  a_2(x_1-x_0)+b_2(\dot{x}_1-\dot{x}_0)\\
\nonumber u_{k} & = & a_2(x_{k+1}-x_k)+b_2(\dot{x}_{k+1}-\dot{x}_k)+a_1(x_{k-1}-x_k)+ \\
\nonumber & & b_1(\dot{x}_{k-1}-\dot{x}_k) \;\;\; \text{ for } \ 1\le k \le N-1\\
 \nonumber u_N & = &  a_1(x_{N-1}-x_N)+b_1(\dot{x}_{N-1}-\dot{x}_N)\, ,
\end{eqnarray}
where $a_1,  a_2,  b_1$ and $b_2$ are constant positive parameters. Plugging \eqref{eq2} into \eqref{system}, we write the dynamics of the configuration error $e_k=x_{k-1}-x_k$ in Laplace domain:
\begin{eqnarray} \label{eq4}
s^2 e_1 & = &  (a_2+b_2s)(e_2-e_1)-(a_1+b_1s)e_1+d'_1\\
\nonumber s^2e_k & = &  (a_2+b_2s)(e_{k+1}-e_k)+(a_1+b_1s)(e_{k-1}-e_k)+d'_k\\
\nonumber && \;\;\; \text{ for } \ 1\le k \le N-1\\
\nonumber s^2e_N & = &  (a_1+b_1s)(e_{N-1}-e_N)-(a_2+b_2s)e_N+d'_N\
  \end{eqnarray}
Here $d'_{k}=d_{k-1}-d_k$ and by the triangle inequality, $\Vert d \Vert_2 < \delta/2$ implies $\Vert d' \Vert_2 < \delta$. A full proof that \eqref{eq4} is stable (before being string stable) has been made, but is left out here due to space constraints. 


\section{Proof of string stability with respect to leader(s)}


\subsection{Analysis I: partial inversion of the dynamics}

The error dynamics \eqref{eq4} can be written compactly as 
\begin{equation} \label{eq6}
SE=D^{'}
\end{equation} 
with matrix $S$ and column vectors $E$, $D'$ given by 
\begin{eqnarray*}
E &=& (e_1, e_2, \dots , e_N)\\
D'&=& (d'_1, d'_2, \dots , d'_N)\\
S &=& \begin{bmatrix}\
    s^2+q & -p_2 & 0 & \dots  & 0 \\
    -p_1 & s^2+q & -p_2 & \dots  & 0 \\
    \vdots & \ddots & \ddots & \ddots & \vdots \\
     0 & \dots & -p_1  & s^2+q & -p_2 \\
    0 & \dots & 0 & -p_1 &    s^2+q
\end{bmatrix}
\end{eqnarray*}
where we have defined the elementary transfer functions $p_1=a_1+b_1s$, $\;\; p_2=a_2+b_2s$ and $q=p_1+p_2$.

%

For a linear system, string stability essentially means: bounded $\Vert D \Vert_2$ implies bounded $\Vert E \Vert_2$, uniformly in $N$. Here we analyze a slightly stronger goal by replacing $D$ with $D'$. This gives a sufficient condition for string stability, as $(L_2,l_2)$-bounded $D$ implies $(L_2,l_2)$-bounded $D'$, uniformly in $N$. When investigating \emph{necessary} conditions for string stability, we will have to restrict the inputs to instances of $(L_2,l_2)$-bounded $D'$ which have a spatial structure that also corresponds to $(L_2,l_2)$-bounded $D$.

To analyze in detail the effect of $D'$ on $E$, we essentially want to invert equation \eqref{eq6}. We will do this in two steps. Namely, first we apply a transformation that makes \eqref{eq6} \emph{almost} diagonal -- i.e.~after transformation each component follows a diagonal dynamics, plus a drive by the boundary vehicles $e_1$ and $e_N$. We are then able to analyze the resulting system by hand. For the first step (transformation), we define the matrix
\begin{eqnarray*}
M &=& \frac{1}{m} \; \begin{bmatrix}
    C & C_2 & C_2^2 & \dots  & C_2^{N-1} \\
    C_1 & C& C_2& \dots  & C_2^{N-2}\\
    \vdots & \ddots & \ddots & \ddots  & \vdots\\
         C_1^{N-2} & \dots & C_1 & C & C_2 \\
    C_1^{N-1} & \dots & C_1^2 & C_1  & C
\end{bmatrix}\\
\text{with } \, m &=& \sqrt{(s^2+q)^2-4p_1p_2} \, ,
\end{eqnarray*}
and $C,C_1,C_2$ to be found. Multiplying both sides of \eqref{eq6} by the proposed matrix $M$, we want to obtain
\begin{equation} \label{eq9}
MSE=QE=MD^{'}
\end{equation} 
with a matrix $Q$ easy to invert. In particular, we impose the structure:
$$Q=MS= \begin{bmatrix}
    q_{1,1} & 0 & 0 & \dots  & q_{1,N} \\
   q_{2,1} & 1 & 0 & \dots  & q_{2,N}\\
    \vdots & \ddots & \ddots & \ddots  & \vdots\\
         q_ {N-1,1}& \dots & 0 & 1 & q_{N-1,N} \\
  q_{N,1} & \dots & 0 & 0  & q_{N,N}\
\end{bmatrix} \, .$$
By working out the matrix multiplication, this imposes the following relations:
\begin{eqnarray}\label{Qcons1}
0 & = & -C_2^k p_2 + C_2^{k+1} (s^2+q) - C_2^{k+2} p_1 \\
\nn 0 & = & -C_1^{k+2} p_2 + C_1^{k+1} (s^2+q) - C_1^k p_1 \\
\nn & & \text{ for } k=1,2,...,N-3 ;\\
\nn 0 & = & -C p_2 + C_2 (s^2+q) - C_2^2 p_1 \\
\nn 0 & = & - C_1^2 p_2 + C_1 (s^2+q) - C p_1 \\
\nn m & = & -C_1 p_2 + C (s^2+q) - C_2 p_1 
\end{eqnarray}
\text{and}
\begin{eqnarray}\label{Qcons2}
m q_{k,N} & = & -C_2^{N-(k+1)} p_2 + C_2^{N-k} (s^2+q)\\
\nn m q_{k+2,1} & = & C_1^{k+1} (s^2+q) - C_1^k p_1 \\
\nn & & \text{ for } k=1,2,...,N-2 ;\\
\nn m q_{1,1} & = & C (s^2+q) - C_2 p_1 \\
\nn m q_{2,1} & = & C_1 (s^2+q) - C p_1 \\
\nn m q_{N-1,N} & = & -C p_2 + (s^2+q) C_2 \\
\nn m q_{N,N} & = & -C_1 p_2 + (s^2+q) C \, .
\end{eqnarray}
The second set of equations \eqref{Qcons2} just defines the $q_{k,1}$ and $q_{k,N}$, to which we will come back later. The first set of equations \eqref{Qcons1} define $C,C_1,C_2$; one checks that they are satisfied if and only if we take
\begin{eqnarray} \label{Cdef}
C & = & 1 \; , \quad C_1 = \frac{(s^2+q) - m}{2 p_2} \; , \quad C_2 = \frac{(s^2+q) - m}{2 p_1} \phantom{kkk} \\
& & \text{with } \;\; m = \sqrt{(s^2+q)^2-4p_1p_2}
\nn  \, .
\end{eqnarray}
In particular, the last line imposes the sign in front of $m$ in the expressions of $C_1$ and $C_2$. To obtain proper transfer functions (\cite{22}), the complex square root of $m$ should be interpreted along the branch for which the dominant $s^2$ terms cancel at high frequencies.

Using \eqref{Qcons1},\eqref{Qcons2} the error dynamics of the vehicles rewrites:
\begin{eqnarray}\label{eq:ErDyns}
e_1 & = & \frac{1}{q_{1,1}} \left(- q_{1,N} e_N + d'_1 / m + \sum_{k=1}^{N-1} C_2^k \, d'_{1+k} / m \right) \\
\nn e_k & = & - q_{k,1} e_1 - q_{k,N} e_N + d'_k / m\\
\nn & & + \sum_{l=1}^{k-1} C_1^l \, d'_{k-l} / m + \sum_{l=1}^{N-k} C_2^l \, d'_{k+l} / m \\
\nn & & \text{ for } k=2,3,...,N-1 \\
\nn e_N & = & \frac{1}{q_{N,N}} \left(- q_{N,1} e_1 + d'_N / m + \sum_{k=1}^{N-1} C_1^k \, d'_{N-k} / m \right) \; .
\end{eqnarray}
We see that the pair $e_1,e_N$ now forms a system of its own, which drives the other vehicles inside the chain. The latter are in addition driven by their local disturbance $d'_k$ and by two flows: a flow of disturbances coming from the front, which we denote
$$f_k = \sum_{l=1}^{k-1} C_1^l \, d'_{k-l} \;\; = \;\; C_1 (f_{k-1}+d'_{k-1}) \; ,$$
and a flow coming from the rear,
$$g_k = \sum_{l=1}^{N-k} C_2^l \, d'_{k+l} \;\; = \;\; C_2 (g_{k+1}+d'_{k+1}) \, .$$

In the next subsection, we analyze separately the parts of $e_k$ related to the disturbance flows and to the $e_1,e_N$ pair. The analysis of the latter brings novel positive news with respect to \cite{MartinecEJC}: no dedicated controller appears to be needed at the boundaries to ensure a well-behaved system.


\subsection{Analysis II: bounding the flow transfer functions}

We first consider the flows $f_k$ and $g_k$. In order to ensure $(L_2,l_2)$ boundedness of those signals, the $H_\infty$ norm of both $C_1$ and $C_2$ would have to be lower than one. We next show that we can tune the controller such that one of those two constraints is satisfied, but not both. We typically choose to have $\Vert C_1(j\omega) \Vert_{\infty} <1$. This leaves the hope of achieving string stability with respect to disturbance inputs e.g.~$d'_1 \neq 0$ on the leading vehicle only. We will then conclude by showing that indeed, assuming $d'_k=0$ for all $k>1$, the $e_1,e_N$ part of the dynamics has an 
$(L_2,l_2)$  bounded influence on the dynamics as well, and thus the asymmetric system can be string stable in that sense.
\vspace{2mm}


\noindent \textbf{Lemma 1:} \emph{Consider the controller \eqref{eq2} with $a_1 \neq 0 \neq a_2$ (no poles cancellation) and $a_1 \neq a_2$. It is impossible to have both $\Vert C_1(j\omega) \Vert_{\infty} \leq 1$ and $\Vert C_2(j\omega) \Vert_{\infty} \leq 1$.}
\vspace{2mm}

\noindent \emph{Proof:} Let us assume $a_2 > a_1$; the converse case is similar. We have
\begin{eqnarray} 
\nonumber C_2 & = & \tfrac{(s^2+p_1+p_2)-\sqrt{(s^2+p_1+p_2)^2 - 4 p_1 p_2}}{2 p_1} \\
\nonumber & = & \tfrac{1}{2p_1} \left( p_1 + s^2 + p_2 - \sqrt{(s^2+p_2-p_1)^2+4 s^2p_1} \right) \\
\nonumber & = & \tfrac{1}{2} \left( 1 + \tfrac{s^2 + p_2}{p_1} - \left(\tfrac{s^2+p_2}{p_1}-1 \right) \sqrt{1 + \tfrac{4 s^2 p_1}{(s^2+p_2-p_1)^2}} \right) \\
\nonumber & \simeq & \tfrac{1}{2} \left( 1 + \tfrac{s^2 + p_2}{p_1} - \left(\tfrac{s^2+p_2}{p_1}-1 \right) - \tfrac{2 s^2}{s^2 + p_2 - p_1} \right) +  O(\vert s \vert^4) \\
\label{eq:C2min} & = & 1 + \tfrac{s^2}{p_1-p_2 - s^2}  +  O(\vert s \vert^4) \, .
\end{eqnarray}
The second line is obtained by square completion. The third line is valid for $\vert s \vert \ll 1$, taking into account that $a_2>a_1$ for the phase of the factor taken out of the square root. The next line is Taylor approximation for the square root for $\vert s \vert \ll 1$; the higher order terms of order $\vert s \vert^4$ can be neglected provided $a_1 \neq a_2$ and $a_1 \neq 0 \neq a_2$, which is the condition to avoid pole cancellation. Replacing $s=j\omega$ in the last line we obtain
$$\vert C_2(j\omega) \vert \simeq \left\vert \tfrac{(a_2-a_1)+(b_2-b_1) j \omega}{(a_2-a_1)+(b_2-b_1) j \omega - \omega^2} \right\vert \; > 1$$
for low frequencies. \hfill $\square$\\

On the positive side, we have the following results.\vspace{2mm}

\noindent \textbf{Lemma 2:} \emph{Consider the controller \eqref{eq2} with $a_1 \neq 0 \neq a_2$ (no poles cancellation).
\newline (a) For any choice of the control parameters we have $\Vert C_1(j\omega) C_2(j\omega) \Vert_\infty \leq 1$.
\newline (b) Taking $p_2 = \alpha p_1$, for any $1\neq\alpha>0$ and any $a_1,b_1>0$, we have $\Vert C_1(j\omega) C_2(j\omega) \Vert_\infty < 1$.
\newline (c) Take case (b) and write $p_1 = \tfrac{\kappa}{1+\alpha} p$ with any $\kappa>0$ and $p=a + b s$, for some fixed $a,b,\kappa>0$. There exists $\bar{\alpha}$ such that for $\alpha>\bar{\alpha}$, we have $\Vert C_1(j\omega) \Vert_{\infty} < 1$.
}
\vspace{2mm}

\noindent \emph{Proof:} (a),(b) We have $|\sqrt{C_1 C_2}| = |1-\sqrt{1-x}| / |\sqrt{x}| =: f(x)$  with $x = \tfrac{4 p_1 p_2}{(s^2+p_1+p_2)^2}$. The property follows from the fact that $f(x) = 1$ for $x \in [1,+\infty)$ and $f(x) < 1$ for all other $x \in \mathbb{C}$. 
For the particular choice of (b), we have $x = \tfrac{4 \alpha p_1^2}{(s^2+(1+\alpha)p_1)^2}$. Since $x(j\omega)$ can be real positive, only if the phases of numerator and denominator match, this can happen only for $(j\omega)^2$ parallel to $(1+\alpha)p_1$, i.e.~$p_1$ real. With $b_1\neq0$ this happens only at $\omega=0$, for which we have $x=4\alpha/(1+\alpha)^2 < 1$. Thus with (b) we never have $x \in [1,+\infty)$, so $f(x) < 1$.

(c) We have
$$C_1 = \frac{s^2 + \kappa p - \sqrt{(s^2+\kappa p)^2-4 \alpha \beta^2 p^2}}{2 \alpha \beta p} \, ,$$
where $\beta = \kappa/(1+\alpha)$. Denote by $g$ the minimum norm of $h(s) = (s^2+\kappa p)^2 / p^2$ over all $s=j\omega$. Recall from standard Bode diagram approximations that $g > 0$ as long as perfect undamped resonance is avoided, i.e.~$b\neq0$. While $h(s)$ stays fixed, we can now decrease the value of $\alpha \beta^2 = \kappa^2 \alpha / (1+\alpha)^2$ to make it arbitrarily smaller than $g$, just by increasing $\alpha$ and decreasing $\beta$ at the same time. This allows to apply the Taylor expansion of $\sqrt{1+x}$ to the square root in $C_1$, uniformly for all $\omega$:
\begin{eqnarray*}
C_1 &=& \frac{\frac{4 \alpha \beta^2 p^2}{(s^2+\kappa p)} + O(\,(\alpha \beta^2/g)^2\,)}{2 \alpha \beta p} \\
& = & \frac{2\kappa}{h(s)(1+\alpha)} + O(\tfrac{\kappa^3 \alpha}{g^2(1+\alpha)^3}).
\end{eqnarray*}
It is clear that the norm of this last expression can be made arbitrarily small by taking $\alpha$ sufficiently large, such that we can make $\Vert C_1(j\omega) \Vert_{\infty}$ smaller than $1$ or in fact than any other value. \hfill $\square$

Lemma 1 indicates that one should not expect $L_2$ string stability with this controller when all the vehicles can be subject to disturbances $d'_k$, except possibly with $a_1=a_2$. This fact has also been established in \cite{Automatica-Knorn} while the present paper was under review. It is not hard to see that, even with other linear controllers having a finite DC gain $a_1=a_2$, it will anyways be impossible to get $(L_2,l_2)$ string stability.


Thanks to Lemma 2(c) however, string stability might hold when disturbances are concentrated on a certain number of leading vehicles, independent of $N$, as is also assumed in \cite{MartinecEJC}. We now further analyze this situation.


\subsection{Analysis III: the $e_1$,$e_N$ subsystem and conclusion}

Let us rewrite the first and last line of \eqref{eq:ErDyns}:
\begin{eqnarray*}
q_{1,1} e_1 & = & -q_{1,N} e_N +  d'_1/m + g_1/m \\
q_{N,N} e_N & = & -q_{N,1} e_1 +  d'_N/m + f_N/m \, .
\end{eqnarray*}
Multiplying the first one by $q_{N,N}$ and substituting the second one into it (respectively conversely), we obtain
\begin{eqnarray*}
\frac{e_1}{d_1^{'}} &=& \frac{mq_{N,N}-mq_{1,N}C_1^{N-1}}{m^2q_{1,1}q_{N,N}-m^2q_{1,N}q_{N,1}} \\
\frac{e_N}{d_1^{'}} &=& \frac{-mq_{N,1}+mq_{1,1}C_1^{N-1}}{m^2q_{1,1}q_{N,N}-m^2q_{1,N}q_{N,1}} \; .
\end{eqnarray*}
provided $q_{N,N}q_{1,1} \neq q_{1,N}q_{N,1}$.
With this expression we can state the following result.
\vspace{2mm}

\noindent \textbf{Theorem 3:} \emph{With appropriate tuning (see Lemma 2), the family of vehicle chains described by the controller \eqref{eq2} for all $N \in \mathbb{N}$, is $(L_2,l_2)$ string stable with respect to disturbances $d'$ \textbf{restricted to the first $\bar{k}$ vehicles only, for some integer $\bar{k}$ independent of $N$}; in other words, it is $(L_2,l_2)$ string stable provided we impose $d'_k=0$ for all $k>\bar{k}$.}
\vspace{2mm}

We will use the following facts later in the proof.
\begin{itemize}
\item[(a)] By choosing $\alpha>\bar{\alpha}$ large enough in the conditions of Lemma 2(c), it is possible to ensure that $m(s)\neq 0$ in the RHP, and thus in particular $m(j\omega)$ bounded away from 0. Indeed, in this setting we have $m = (s^2+p) \sqrt{1-\frac{4\alpha}{(1+\alpha)^2}\,(\frac{p}{s^2+p})^2}$. The second-order polynomial $s^2+p$ has all roots in LHP for positive coefficients. Since moreover $\frac{p}{s^2+p}$ goes to $0$ for $|s|$ going to infinity, we can upper bound $|\frac{p}{s^2+p}|^2<\bar{\eta}$ in the RHP. Then by taking $\alpha$ large enough, we can make $\frac{4\alpha}{(1+\alpha)^2}$ small enough, in particular such that $\sqrt{1-\frac{4\alpha}{(1+\alpha)^2} \bar{\eta}} > 0$, thus implying the property.
\item[(b)] For a tuning as in \emph{Lemma 2(c)}, we can give a lower bound $\eta_2>0$ for the norm of $(\tfrac{s^2+q+m}{2})^2$ in the RHP. Indeed, note that
$s^2+q+m = (s^2+p) + (s^2+p) \sqrt{1-\frac{4\alpha}{(1+\alpha)^2}\,(\frac{p}{s^2+p})^2}$. The factor $s^2+p$ has roots in LHP, like for point (a). We have also explained in point (a) that by choosing $\alpha$ large enough, we can make the term $\frac{4\alpha}{(1+\alpha)^2}\,(\frac{p}{s^2+p})^2$ arbitrarily small in the RHP. It is then clear that we can ensure $1+\sqrt{1-\frac{4\alpha}{(1+\alpha)^2}\,(\frac{p}{s^2+p})^2} \neq 0$ in the RHP.
\end{itemize}

\noindent \emph{Proof:} The basic case is of course when $\bar{k}=1$ i.e.~only the leader is subject to a disturbance. We here provide the proof for this case; the general case is similar.

We will choose $p_1,p_2$ according to \emph{Lemma 4.4(c)} such that $\Vert C_1(j\omega) C_2(j\omega) \Vert_\infty < 1$, and with the controller parameterized via $\alpha$ and $p$.

We thus assume $d'_k = 0$ for all $k>1$, which implies $g_k = 0$ for all $k$ and $f_k = C_1^{k-1} \, d'_1$. 

We first analyze $e_1$. A few computations lead to
\begin{eqnarray*}
\frac{e_1}{d'_1}=\frac{1-(C_1C_2)^N}{\frac{s^2+q+m}{2}[1-{(C_1C_2)}^{N+1}]}.
& =: & H_1(s) =: G_1(s) \, .
\end{eqnarray*}
A first point is to prove stability of $H_1(s)$. By the property (b) above, this comes down to proving that 
$(C_1 C_2)^{N+1} \neq 1$ in the RHP. Since $C_1 C_2 = (\frac{s^2+q-m}{2})\, /\, (\frac{s^2+q+m}{2})$, we have to prove that
$$(\frac{s^2+q+m}{2})^N \neq (\frac{s^2+q-m}{2})^N \, .$$
As $N$ can take arbitrary integer values, we will show that $\vert s^2+q+m \vert \neq \vert s^2+q-m \vert$. To have equality, we would need $m$ perpendicular to $s^2+q$ in the complex plane. But analyzing $m$ as in property (a) above, we can choose $\alpha$ such that $m = (s^2 + p) \sqrt{1 + \eta_3}$ with $|\eta_3| \ll 1$, such that perpendicularity cannot be achieved. Thus, $H_1$ is stable.

We next check string stability. For large $\omega$ we have $\Vert C_1(j\omega) C_2(j\omega) \Vert = O(1/\omega^2)$, so $H_1(s)$ behaves like $ O(\omega^2) \, / \, O(\omega^4)$ for large $\omega,N$, with leading coefficients independent of $N$. For any $\xi>0$, we can thus define $\bar{\omega}$ such that $\vert H_1(j\omega) \vert < \xi$ for all $\omega > \bar{\omega}$ and for all $N>3$. For the compact domain $\omega<\bar{\omega}$, thanks to property (b) above and to $\Vert C_1(j\omega) C_2(j\omega) \Vert<1$, we have a bound on $\Vert H_1(j\omega) \Vert_\infty$ which is independent of $N$.\vspace{2mm}

We next turn to $e_N$.
\newline Similarly we have
\begin{eqnarray*}
\frac{e_N}{d_1^{'}}&=&\frac{(\frac{s^2+q+m}{2})C_1^{N-1}+[C_1^{N-2}p_1-C_1^{N-1}(s^2+q)]}{(\frac{s2+q+m}{2})^2-(C_1C_2){N-2}p_1p_2[1-\frac{(s^2+q-m)(s^2+q)}{2p_1p_2}]}\\
&=&\frac{(\frac{s^2+q+m}{2})C_1^{N-1}+p_1C_1^{N-2}[1-\frac{(s^2+q-m)(s^2+q)}{2p_1p_2}]}{(\frac{s^2+q+m}{2})^2-(C_1C_2)^{N-2}p_1p_2[1-\frac{(s^2+q-m)(s^2+q)}{2p_1p_2}]}\\
&=&\frac{mC_1C_2}{p_2(\frac{s^2+q+m}{2})}\cdot \frac{C_1^{N-2}}{1-(C_1C_2)^{N+1}}=:H_N(s)\; .
\end{eqnarray*}
By the same arguments the transfer function $H_N$ is stable and the transfer function $G_N := H_N / C_1^{N-2}$ from $C_1^{N-2} d'_1$ to $e_N$ is bounded independently of $N$.\vspace{2mm}

For the other vehicles, we then have
\begin{eqnarray*}
\frac{e_k}{d_1^{'}}&=&-q_{k,1}\frac{e_1}{d_1^{'}}-q_{k,N}\frac{e_N}{d_1^{'}}+C_1^{k-1}/m\\
&=&C_1^{k-2}\left(\frac{C_1}{m}-\frac{[1-(C_1C_2)^N]C_1C_2(\frac{s^2+q-m}{2})+m(C_1C_2)^{N-k+2}}{mp_2[1-(C_1C_2)^{N+1}]}\right)\\
 &&=: H_k(s) =: G_k(s) C_1^{k-2} \; .
\end{eqnarray*}
Proving stability involves the same elements as for vehicle $1$, plus requiring $m\neq 0$ in the RHP; the latter property is proved in item (a) above. Towards proving string stability, one can also apply the same arguments as for $G_1(s)$ to the different terms of $G_k(s)$: they are bounded for $\omega \gg 1$, and for finite $\omega$ we can bound $\Vert G_k(j\omega)\Vert_\infty$, independently of $N$ and of $k$, provided we have a lower bound on $\Vert m(j\omega) \Vert_\infty$. The latter is also ensured by property (a) above.

Taking all things together, we have
\begin{eqnarray*}
\Vert e(\cdot) \Vert_2^2 & \leq & \sum_{k=1}^N  \Vert H_k(j\omega) \Vert_{\infty}^2 \; \Vert d'_1(\cdot) \Vert_2^2 \\
&= & \sum_{k=2}^N  \Vert G_k(j\omega) \Vert_{\infty}^2 \; \Vert C_1(j\omega)^{k-2} \Vert_{\infty}^2 \; \Vert d'_1(\cdot) \Vert_2^2 \\
&& + \Vert G_1(j\omega) \Vert_{\infty}^2 \; \Vert d'_1(\cdot) \Vert_2^2 \\
& \leq & \Vert d'_1(\cdot) \Vert_2^2 \; \Vert G_{\max}(j\omega) \Vert_{\infty}^2 \; r_N \; .
\end{eqnarray*}
Here $G_{\max}$ is the transfer function, among the $G_k$, with the largest $H_\infty$ norm; we have just shown that this norm is bounded independently of $N$. And
$$r_N := 1 + \sum_{k=2}^{N} r^{2(k-2)} = 1 + \frac{1-r^{2(N-1)}}{1-r^2}$$
with $r:=\Vert C_1(j\omega) \Vert_{\infty}$ is bounded independently of $N$ when $r<1$; the latter condition can be satisfied by \emph{Lemma 2(c)}. This concludes the proof. \hfill $\square$


\section{Simulations}\label{4}

\begin{figure}[h!]
\includegraphics[width=0.49\textwidth, trim=12mm 0mm 10mm 7mm, clip=true]{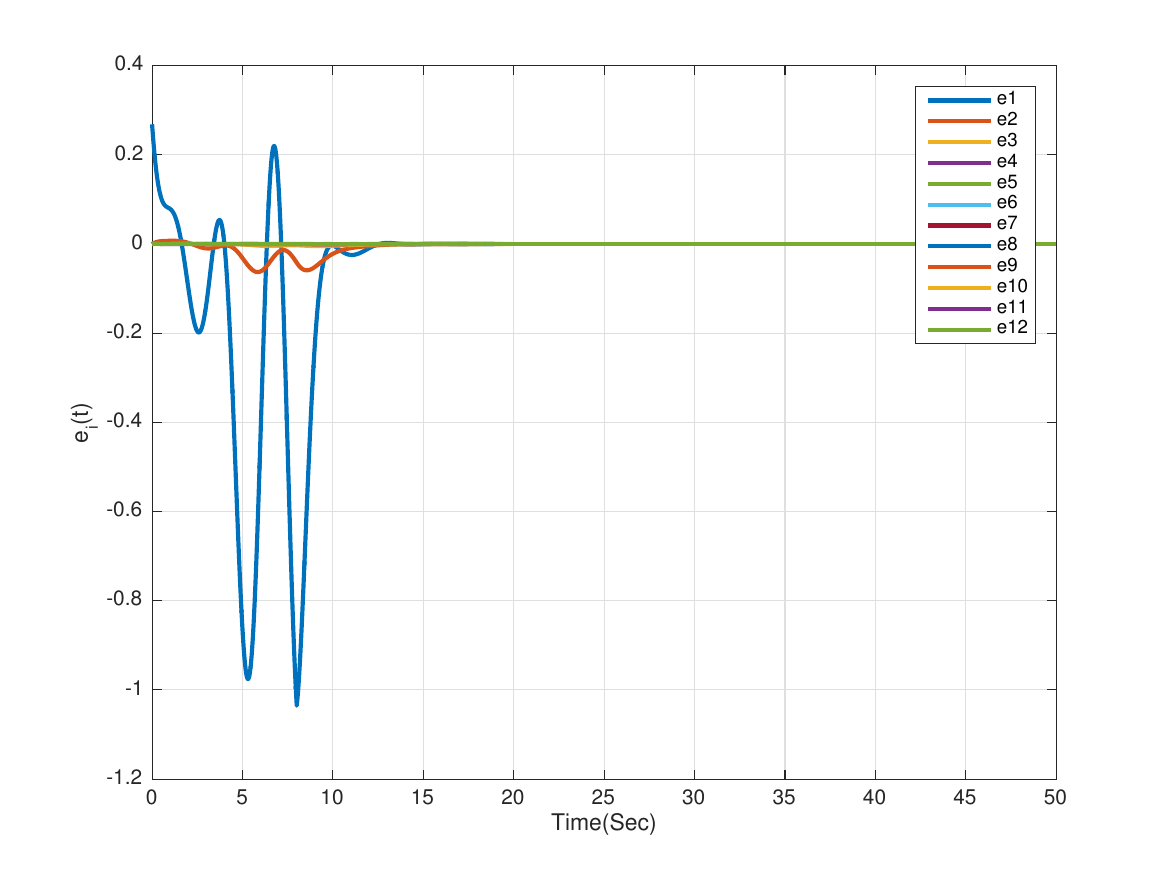}  
\caption{Spacing errors $e_i(t)$ of a platoon with 12 following vehicles.}\label{fig3}
\end{figure}

\begin{figure}[h!]
\includegraphics[width=0.49\textwidth, trim=12mm 0mm 10mm 7mm, clip=true]{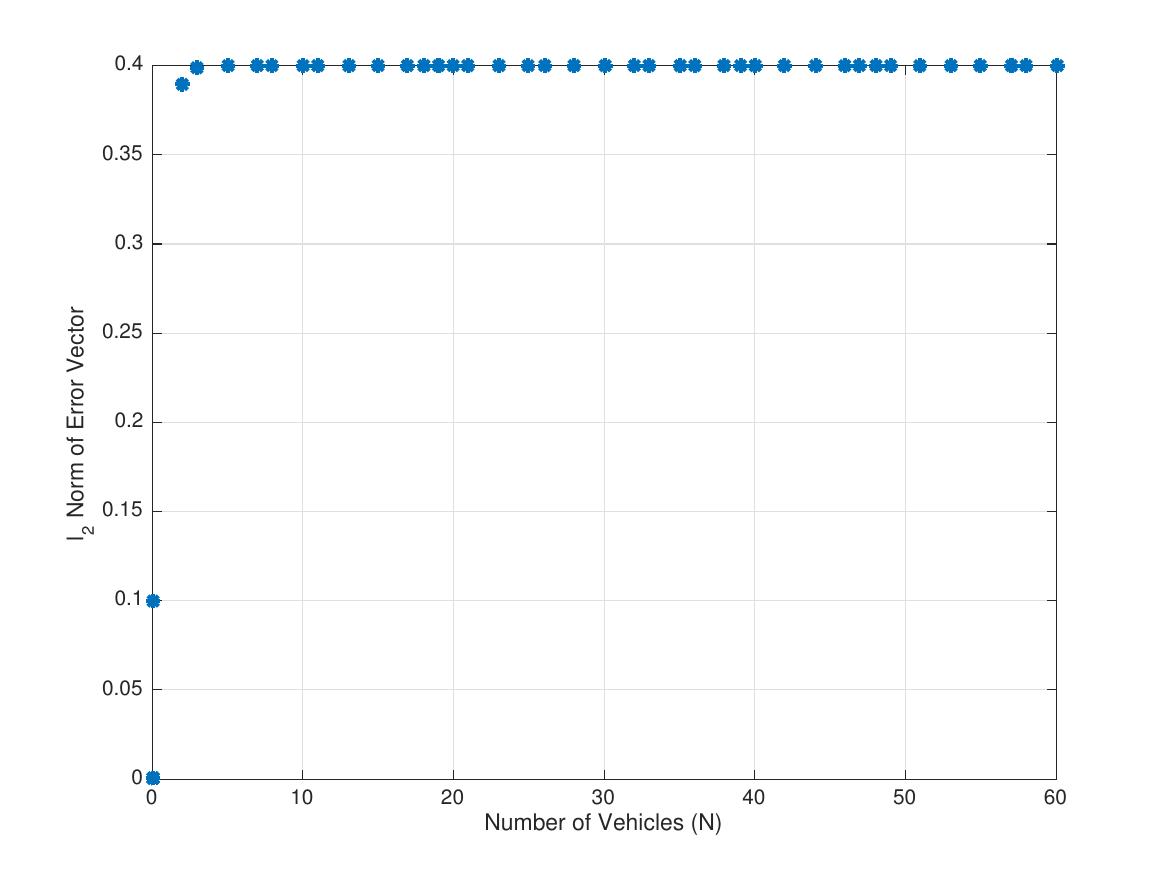}  
\caption{String stability criterion as a function of chain length $N$: the $(L_2,l_2)$ norm of $(e_1(t),...,e_N(t))$ indeed stays bounded as a function of $N$.}\label{fig4}
\end{figure}

We can briefly illustrate the effectiveness of the proposed asymmetric bidirectional controller in simulation. We apply a short disturbance on the leading vehicle of a platoon with control parameters $a_1=1$, $b_1=1$, $a_2=10$ and $b_2=100$. This is not exactly the ``practical'' tuning $p_2 = \alpha \, p_1$ exploited in the proof, but it appears to work as well, showing some (expected) robustness with respect to the tuning parameters. Figure \ref{fig3} shows the evolution in time of the spacing errors $e_i(t)$, for a network of 12 vehicles. It is apparent that the error decreases not only in time but also along the vehicle chain -- after 3 vehicles essentially, it becomes barely visible.
Figure \ref{fig4} confirms that this controller satisfies the definition of string stability, by showing that the $(L_2,l_2)$-norm of the error vector, as a function of the length $N$ of the chain, converges to a constant bound.

\section{Conclusion}\label{5}

In this paper, we have shown that introducing asymmetry in bidirectional controllers can provide concrete benefits also towards $(L_2,l_2)$ string stability. More precisely, we have shown that a simple asymmetric coupling among vehicles allows to solve this string stability problem for a vehicle chain of length $N$, provided the disturbances are acting on a few ($N$-independent) leading vehicles only. We have also re-proved, with this alternative flow formulation, that if disturbances act on all vehicles with such controller, then no parameter values can achieve string stability. A straightforward extension satisfying string stability would be to allow disturbances $d'_k$ that decrease exponentially with $k$, at the same rate as the $C_2$ function in our analysis increases. Future work will concentrate on finding minimal alternatives to the present setting, possibly exploiting the property $\Vert C_1 C_2 \Vert_{\infty(AC)} \leq 1$, in order to solve string stability under arbitrary disturbances.


\end{document}